\documentstyle[12pt]{article}
\begin{document}
\newpage
\pagestyle{empty}
\setcounter{page}{1}
%%%%%%%%%%%%%%%%%%%%%%%%%%%%%%%
%%%%%%%%%%%%%%%%%%%%%%%%%%%%%%
%%%%%%%%%%%%%%%%%%%%%%%
%\renewcommand{\thesection}{\Roman{section}}
%\renewcommand{\theequation}{\thesection.\arabic{equation}}
%\newcommand{\sect}[1]{\setcounter{equation}{0}\section{#1}}
%%%%%%%%%%%%%%%%%%%%%%%%
%%%%%%%%%%%%%%%%%%%%%%%%%%%
%%%%%%%%%%%%%%%%%%%%%%%%%%%%%%
%%% ** start of amsfont definitions **
\newfont{\twelvemsb}{msbm10 scaled\magstep1}
\newfont{\eightmsb}{msbm8} \newfont{\sixmsb}{msbm6} \newfam\msbfam
\textfont\msbfam=\twelvemsb \scriptfont\msbfam=\eightmsb
\scriptscriptfont\msbfam=\sixmsb \catcode`\@=11
\def\Bbb{\ifmmode\let\next\Bbb@\else \def\next{\errmessage{Use
      \string\Bbb\space only in math mode}}\fi\next}
\def\Bbb@#1{{\Bbb@@{#1}}} \def\Bbb@@#1{\fam\msbfam#1}
\newfont{\twelvegoth}{eufm10 scaled\magstep1}
\newfont{\tengoth}{eufm10} \newfont{\eightgoth}{eufm8}
\newfont{\sixgoth}{eufm6} \newfam\gothfam
\textfont\gothfam=\twelvegoth \scriptfont\gothfam=\eightgoth
\scriptscriptfont\gothfam=\sixgoth \def\frak{\frak@}
\def\frak@#1{{\fam\gothfam{{#1}}}} \def\frak@@#1{\fam\gothfam#1}
\catcode`@=12
%%% ** end of amsfont definitions **
%
%
%
%\def\Bbb{\bf}
\def\CC{{\Bbb C}}
\def\NN{{\Bbb N}}
\def\QQ{{\Bbb Q}}
\def\RR{{\Bbb R}}
\def\ZZ{{\Bbb Z}}
\def\cA{{\cal A}}          \def\cB{{\cal B}}          \def\cC{{\cal C}}
\def\cD{{\cal D}}          \def\cE{{\cal E}}          \def\cF{{\cal F}}
\def\cG{{\cal G}}          \def\cH{{\cal H}}          \def\cI{{\cal I}}
\def\cJ{{\cal J}}          \def\cK{{\cal K}}          \def\cL{{\cal L}} 
\def\cM{{\cal M}}          \def\cN{{\cal N}}          \def\cO{{\cal O}}
\def\cP{{\cal P}}          \def\cQ{{\cal Q}}          \def\cR{{\cal R}} 
\def\cS{{\cal S}}          \def\cT{{\cal T}}          \def\cU{{\cal U}}
\def\cV{{\cal V}}          \def\cW{{\cal W}}          \def\cX{{\cal X}}
\def\cY{{\cal Y}}          \def\cZ{{\cal Z}}
\def\qed{\hfill \rule{5pt}{5pt}}
\def\arcsinh{\mathop{\rm arcsinh}\nolimits}
\newtheorem{theorem}{Theorem}
\newtheorem{prop}{Proposition}
\newtheorem{conj}{Conjecture}
\newenvironment{result}{\vspace{.2cm} \em}{\vspace{.2cm}}

\pagestyle{plain}

\begin{center}

  {\Large {\bf JORDANIAN QUANTUM ALGEBRA ${\cal U}_{\sf h}(sl(N))$ 
VIA CONTRACTION METHOD AND MAPPING}}
\\[1cm]

{B. ABDESSELAM$^{\dagger,}$\footnote{Permanent address: 
Laboratoire de Physique Th\'eorique, Centre Universitaire Mustapha 
Stambouli, 29000-Mascara, Alg\'erie.}$^{,}$\footnote{E-mail: 
boucif@celfi.phys.univ-tours.fr}, 
A. CHAKRABARTI$^{\ddagger,}$\footnote{E-mail: chakra@cpth.polytechnique.fr} 
and R. CHAKRABARTI$^{*,}$\footnote{Permanent address: Department of 
Theoretical Physics, University of Madras, Guindy Campus, 
Madras 600025, India.}$^{,}$\footnote{E-mail: 
ranabir@imsc.ernet.in}}

\smallskip  

\textit{$^{\dagger}$
Laboratoire de Math\'ematique et Physique Th\'eorique, 
Facult\'e des Sciences et Techniques, \\ 
Parc de Grandmont-F 37041 Tours, France.\\
$^{\ddagger}$Centre de 
Physique Th\'eorique, Ecole Polytechnique, 
91128 Palaiseau Cedex, France.\\
$^{*}$Institute of Mathematical Sciences, 
Madras 600113, India. }

\end{center}
\smallskip
\smallskip

\centerline{May 2001}

\begin{abstract}

\noindent{}Using the contraction procedure introduced by us in Ref. 
\cite{ACC2}, we construct, in the first part of the present letter, the 
Jordanian quantum Hopf algebra ${\cal U}_{\sf h}(sl(3))$ which has  
a remarkably simple coalgebraic structure and contains the Jordanian Hopf 
algebra ${\cal U}_{\sf h}(sl(2))$, obtained by Ohn, as a subalgebra. A 
nonlinear map between ${\cal U}_{\sf h}(sl(3))$ and the classical 
$sl(3)$ algebra is then established. In the second part, we give the 
higher dimensional Jordanian algebras ${\cal U}_{\sf h}(sl(N))$ for 
all $N$. The Universal ${\cal R}_{\sf h}$-matrix of ${\cal U}_{\sf h}
(sl(N))$ is also given.

\smallskip

\noindent{\bf Keywords:} Standard quantization, Nonstandard quantization, 
contraction procedure, Hopf algebra, universal ${\cal R}$-matrix, 
Irreducible representations (irreps.).

\end{abstract}

%\tableofcontents

%\newpage 

\section{Introduction}

It is well known that the enveloping Lie algebra ${\cal U}(sl(N))$ has two 
quantizations: The first one called the {\it Drinfeld-Jimbo deformation} 
or the {\it standard quantum deformation} \cite{D1,J} is quasitriangular 
(${\cal R}_{21}{\cal R}\neq I$), whereas the second one called the {\it
Jordanian deformation} or the {\it non-standard quantum deformation} 
\cite{DMMZ} is triangular (${\cal R}_{21}{\cal R}= I$). A typical example of 
Jordanian quantum algebras was first introduced by Ohn \cite{O}. In general, 
nonstandard quantum algebras are obtained by applying Drinfeld twist to the 
corresponding Lie algebras \cite{D2}. The twisting that produces an algebra 
isomorphic  to the Ohn algebra \cite{O} is found in \cite{Og,GGS}. 

Recently, the twisting procedure was extensively employed  to study a wide 
variety of Jordanian deformed algebras, such as  
${\cal U}_{\sf h}(sl(N))$ algebras 
\cite{KLM,LDO1,LDO2,LMDO}, symplectic algebras ${\cal U}_{\sf h}(sp(N))$ 
\cite{AKL}, orthogonal algebras ${\cal U}_{\sf h}(so(N))$ 
\cite{KLDO,KL,KLS,LSK} and orthosymplectic superalgebra ${\cal U}_{\sf h}
(osp(1|2))$ \cite{CK,K}. It appears from these studies that:

\smallskip

{\bf 1.} The non-standard quantum algebras have undeformed commutation 
relations;

\smallskip

{\bf 2.} The Jordanian deformation appear only in the coalgebraic 
structure;

\smallskip

{\bf 3.} The coproduct and the antipode maps have very complicated 
forms in comparison with the Drinfeld-Jimbo and the Ohn deformations. 

\smallskip

To our knowledge, Jordanian quantum algebra ${\cal U}_{\sf h}(sl(N))$ has 
been written explicitly, with a simple coalgebra, only for $N=2$ \cite{O}. 
The main object of the present letter is to construct the Jordanian quantum 
algebra ${\cal U}_{\sf h}(sl(3))$ using the contraction procedure developed 
in \cite{ACC2} and the map studied in Refs. \cite{ACC2,ACCS}. The 
${\cal U}_{\sf h}(sl(3))$ algebra presented here has the following properties:

\smallskip

{\bf 1.} The Ohn algebra ${\cal U}_{\sf h}(sl(2))$ is included in our 
structure ${\cal U}_{\sf h}(sl(3))$ in a natural way as a Hopf
subalgebra and appear 
here from the longest root generators {\it i.e.} from $e_3$, $f_3$ 
and their corresponding Cartan generator $h_3$;

\smallskip

{\bf 2.} Our Jordanian deformed ${\cal U}_{\sf h}(sl(3))$ algebra may
be regarded as the dual Hopf algebra of the function algebra
$Fun_{\sf h}(SL(3))$ studied in \cite{Ali}; 

\smallskip

{\bf 3.} The present ${\cal U}_{\sf h}(sl(3))$ algebra is 
endowed with a relatively simple coalgebra structure (as compared to 
previous studies \cite{KLM,LDO1,LDO2,LMDO}).

\smallskip

Implementing our contraction technique we subsequently obtain higher
dimensional Jordanian quantum algebras ${\cal U}_{\sf h}(sl(N))$ for 
arbitrary values of $N$.  

This letter is organized as follows: The Jordanian quantum algebra 
${\cal U}_{\sf h}(sl(3))$ is introduced via a nonlinear map and proved 
to be a Hopf algebra in section 2. The irreducible representations (irreps.) 
of  ${\cal U}_{\sf h}(sl(3))$ are also given. Higher dimensional algebras 
${\cal U}_{\sf h}(sl(N))$, $N\geq 4$ are presented in the sections 3 and 4.

\section{${\cal U}_{\sf h}(sl(3))$: Map, Hopf Algebra, Irreps. and 
${\cal R}_{\sf h}$-matrix} 

In this letter, ${\sf h}$ is an arbitrary complex number. It was proved 
in \cite{ACC2} that the ${\cal R}_{\sf h}$-matrix of the Jordanian 
quantum algebra ${\cal U}_{\sf h}(sl(3))$ can be obtained from the ${\cal 
R}_q$-matrix associated to the Drinfeld-Jimbo quantum algebra 
${\cal U}_q(sl(3))$ through a specific contraction which is singular in the  
$q\rightarrow 1$ limit. For the transformed matrix, the singularities, 
however, cancel yielding a well-defined construction.  Here we 
assume the ${\cal U}_q(sl(3))$ Hopf algebra to be well-known \cite{M}. 

For brevity and simplicity we limit ourselves to 
(fundamental irrep.)$\;\otimes\;$(arbitrary irrep.). Recall that for
${\cal U}_q(sl(3))$ algebra the $R_q$-matrix in the representation 
(fund.)$\;\otimes\;$(arb.) reads \cite{M}: 
\begin{eqnarray}
&& R_q=\biggl(\pi_{(fund.)}\;\otimes\;\pi_{(arb.)}\biggr)
{\cal R}_q\nonumber \\
&& \phantom{R_q}=
\pmatrix{q^{{1\over 3}(2h_1+h_2)} 
& q^{{1\over 3}(2h_1+h_2)}\Lambda_{12}  
& q^{{1\over 3}(2h_1+h_2)}\Lambda_{13}
\cr
0 
& q^{-{1\over 3}(h_1-h_2)} 
& q^{-{1\over 3}(h_1-h_2)}\Lambda_{23}  
\cr
0 
& 0 
& q^{-{1\over 3}(h_1+2h_2)}  \cr}, 
\end{eqnarray}
where 
\begin{eqnarray}
&&\Lambda_{12}=q^{-1/2}\bigl(q-q^{-1}\bigr)q^{-h_1/2}{\hat f}_1,
\nonumber\\
&&\Lambda_{13}=q^{-1/2}\bigl(q-q^{-1}\bigr){\hat f}_3
q^{-{1\over 2}(h_1+h_2)},
\nonumber\\
&&\Lambda_{23}=q^{-1/2}\bigl(q-q^{-1}\bigr)q^{-h_2/2}{\hat f}_2. 
\end{eqnarray}
The elements $k_1^{\pm 1}=q^{\pm h_1}$, $k_2^{\pm 1}=q^{\pm h_2}$, 
$k_3^{\pm 1}=q^{\pm h_3}
=q^{\pm (h_1+ h_2)}$, ${\hat e}_1$, ${\hat e}_2$, ${\hat e}_3={\hat e}_1
{\hat e}_2-q^{-1}{\hat e}_2{\hat e}_1$, ${\hat f}_1$, 
${\hat f}_2$ and ${\hat f}_3={\hat f}_2{\hat f}_1-q{\hat f}_1{\hat f}_2$ 
are the ${\cal U}_q(sl(3))$ generators. The corresponding classical 
generators are denoted by 
$h_1$, $h_2$, $h_3=h_1+h_2$, $e_1$, $e_2$, $e_3=e_1e_2-e_2e_1$, $f_1$, $f_2$ 
and $f_3=f_2f_1-f_1f_2$. 

We have shown in \cite{ACC2} that the nonstandard 
$R_{\sf h}$-matrix (in the representation (fund.)$\;\otimes\;$(arb.)) arise 
from the $R_q$-matrix (in (fund.)$\;\otimes\;$(arb.)) 
as follows:  
\begin{eqnarray}
&& R_{\sf h}=\lim_{q\rightarrow 1}
\biggl[E_q\biggl({{\sf h}{\hat e}_3\over q-1}\biggr)_{(fund.)}\otimes 
E_q\biggl({{\sf h}{\hat e}_3\over q-1}\biggr)_{(arb.)}\biggr]^{-1}
R_q\biggl[E_q\biggl({{\sf h}{\hat e}_3\over q-1}\biggr)_{(fund.)}
\otimes E_q\biggl({{\sf h}{\hat e}_3\over q-1}\biggr)_{(arb.)}\biggr]
\nonumber\\
&&\phantom{R_{\sf h}}=\lim_{q\rightarrow 1}\pmatrix{
E_q^{-1}\bigl({{\sf h}{\hat e}_3\over q-1}\bigr) & 0 & 
-{{\sf h}\over q-1}E_q^{-1}\bigl({{\sf h}{\hat e}_3\over q-1}\bigr) \cr
0 & E_q^{-1}\bigl({{\sf h}{\hat e}_3\over q-1}\bigr)& 0\cr
0 & 0 & E_q^{-1}\bigl({{\sf h}{\hat e}_3\over q-1}\bigr)\cr}
R_q \pmatrix{
E_q\bigl({{\sf h}{\hat e}_3\over q-1}\bigr) & 0 & 
{{\sf h}\over q-1}E_q\bigl({{\sf h}{\hat e}_3\over q-1}\bigr) \cr
0 & E_q\bigl({{\sf h}{\hat e}_3\over q-1}\bigr)& 0\cr
0 & 0 & E_q\bigl({{\sf h}{\hat e}_3\over q-1}\bigr)\cr} \nonumber \\
&& \phantom{R_{\sf h}}=\pmatrix{T 
&2{\sf h}T^{-1/2}e_2 & -{{\sf h}\over 2}(T+T^{-1})\bigl(h_1+h_2\bigr)
+{{\sf h}\over 2}\bigl(T-T^{-1}\bigr)\cr
0 & I & -2{\sf h}T^{1/2}e_1  \cr
0 & 0 & T^{-1}  \cr},
\end{eqnarray} 
where 
\begin{eqnarray}
&& T={\sf h}e_3+\sqrt{1+{\sf h}^2e_3^2}, \qquad\qquad
T^{-1}=-{\sf h}e_3+\sqrt{1+{\sf h}^2e_3^2} .
\end{eqnarray}
The deformed exponential in (3) is defined by
\begin{eqnarray}
&& E_q(x)=\sum_{n=0}^{\infty}{x^n\over [n]!}, \nonumber \\
&& [n]={q^n-q^{-n}\over q-q^{-1}},\qquad\qquad [n]!=[n]\times [n-1]!,
\qquad\qquad
[0]!=1. 
\end{eqnarray}

The following properties can be pointed out:

\smallskip

{\bf 1.} The corner elements of (3) have exactly the same structure as in the
$R_{\sf{h}}$-matrix of ${\cal U}_{\sf{h}}(sl(2))$. This implies that the 
classical generators $e_3$, $h_3=h_1+h_2$ and $f_3$ of ${\cal U}(sl(3))$ are 
deformed (for the nonstandard quantization: ${\cal U}(sl(3))\longrightarrow
{\cal U}_{\sf{h}}(sl(3))$) as follows \cite{ACC2,ACCS}:
\begin{eqnarray}
&& T={\sf{h}}e_3+\sqrt{1+{\sf{h}}^2e_3^2}, \qquad\qquad\;\; 
T^{-1}=-{\sf{h}}e_3+\sqrt{1+{\sf{h}}^2e_3^2}, \nonumber \\
&& H_3=\sqrt{1+{\sf{h}}^2e_3^2}h_3, \qquad \qquad \qquad
F_3=f_3-{{\sf{h}}^2\over 4}e_3\bigl(h_3^2-1\bigr), 
\end{eqnarray}
and evidently satisfy the commutation relations \cite{O}
\begin{eqnarray}
&& TT^{-1}=T^{-1}T=1, \nonumber \\
&& [H_3,T]=T^2-1, \qquad \qquad \qquad \;[H_3,T^{-1}]=T^{-2}-1, \nonumber \\
&& [T,F_3]={{\sf{h}}\over 2}\biggl(H_3T+TH_3\biggr),\qquad
 [T^{-1},F_3]=-{{\sf{h}}\over 2}\biggl(H_3T^{-1}+T^{-1}H_3\biggr), \nonumber\\
&& [H_3,F_3]=-{1\over 2}\biggl(TF_3+F_3T+T^{-1}F_3
+F_3T^{-1}\biggr).
\end{eqnarray}
With the following definition (see Ref. \cite{O})
\begin{eqnarray}
E_3={\sf h}^{-1}\ln T={\sf h}^{-1}\arcsinh{\sf h}e_3,
\end{eqnarray}
it follows that the elements $H_3$, $E_3$ and $F_3$ satisfy the relations
\begin{eqnarray}
&& [H_3,E_3]= 2 {\sinh {\sf h}E_3\over {\sf h}}, \nonumber\\
&& [H_3,F_3]=-F_3\biggl(\cosh {\sf h}E_3\biggr)-
\biggl(\cosh {\sf h}E_3\biggr)F_3,\nonumber \\  
&& [E_3,F_3]=H_3,
\end{eqnarray}
where it is obvious that as ${\sf h}\longrightarrow 0$, we have 
$(H_3,E_3,F_3)\longrightarrow (h_3,e_3,f_3)$. 
It is now evident from (7) that ${\cal U}_{\sf{h}}(sl(2))\subset 
{\cal U}_{\sf{h}}(sl(3))$.

\smallskip

{\bf 2.} The expression (3) of the $R_{\sf h}$-matrix indicates that the 
simple root generators $e_1$ and $e_2$ are deformed as follows:
\begin{eqnarray}
&& E_1=\sqrt{{\sf h}e_3+\sqrt{1+{\sf h}^2e_3^2}}e_1=T^{1/2}e_1,\nonumber \\
&& E_2=\sqrt{{\sf h}e_3+\sqrt{1+{\sf h}^2e_3^2}}e_2=T^{1/2}e_2.
\end{eqnarray}
To complete our ${\cal U}_{\sf h}(sl(3))$ algebra, we introduce the 
following ${\sf h}$-deformed generators:
\begin{eqnarray}
&& F_1=\sqrt{-{\sf h}e_3+\sqrt{1+{\sf h}^2e_3^2}}f_1
       +{{\sf h}\over 2}\sqrt{{\sf h}e_3
       +\sqrt{1+{\sf h}^2e_3^2}}e_2h_3=
       T^{-1/2}\biggl(f_1+{{\sf h}\over 2}e_2Th_3\biggr), \nonumber \\
&& F_2=\sqrt{-{\sf h}e_3+\sqrt{1+{\sf h}^2e_3^2}}f_2
       -{{\sf h}\over 2}\sqrt{{\sf h}e_3
       +\sqrt{1+{\sf h}^2e_3^2}}e_1h_3=
       T^{-1/2}\biggl(f_2-{{\sf h}\over 2}e_1Th_3\biggr), \nonumber \\
&& H_1=\biggl(-{\sf h}e_3+\sqrt{1+{\sf h}^2e_3^2}\biggr)
       \biggl(\sqrt{1+{\sf h}^2e_3^2}h_1
       +{{\sf h}\over 2}e_3(h_1-h_2)\biggr)=
       h_1-{{\sf h}\over 2}e_3T^{-1}h_3, \nonumber\\
&& H_2=\biggl(-{\sf h}e_3+\sqrt{1+{\sf h}^2e_3^2}\biggr)
       \biggl(\sqrt{1+{\sf h}^2e_3^2}h_2
       -{{\sf h}\over 2}e_3(h_1-h_2)\biggr)=
       h_2-{{\sf h}\over 2}e_3T^{-1}h_3.
\end{eqnarray}
The expressions (6), (10) and (11) constitute a realization of the Jordanian 
algebra ${\cal U}_{\sf h}(sl(3))$ with the classical generators via a nonlinear
map. This immediately yields the irreducible representations (irreps.) of 
${\cal U}_{\sf h}(sl(3))$ in an explicit and simple manner. 

\begin{prop}
The Jordanian 
algebra ${\cal U}_{\sf h}(sl(3))$ is an associative algebra over $\CC$ 
generated by $H_1$, $H_2$, $H_3$, $E_1$, $E_2$, $T$, $T^{-1}$, $F_1$, $F_2$ 
and $F_3$, satisfying, along with (7), the commutation relations
\begin{eqnarray}
&& [H_1,H_2]=0, \qquad\qquad\qquad\qquad\qquad\qquad\qquad
     [H_1,T^{-1}H_3]=[H_2,T^{-1}H_3]=0,\nonumber \\
&& [H_1,E_1] =2E_1,\qquad\qquad\qquad\qquad\qquad\qquad\;\;\;
     [H_2,E_2] =2E_2,\nonumber \\
&& [H_1,E_2]=-E_2,\qquad\qquad\qquad\qquad\qquad\qquad\;\; 
     [H_2,E_1]=-E_1, \nonumber\\ 
&& [T^{-1}H_3,E_1]=E_1,\qquad\qquad\qquad\qquad\qquad\;\;\;\;\;\;
     [T^{-1}H_3,E_2]=E_2, \nonumber\\
&& [H_1,F_1] =-2F_1+{\sf h}E_2T^{-1}H_3, \qquad\qquad\qquad\;
     [H_2,F_2] =-2F_2-{\sf h}E_1T^{-1}H_3, \nonumber\\
&& [H_1,F_2] =F_2-{\sf h}E_1T^{-1}H_3, \qquad\qquad\qquad \;\;\;\;\;
     [H_2,F_1] =F_1+{\sf h}E_2T^{-1}H_3, \nonumber\\
&& [TH_3,F_1] =-T^2F_1, \qquad\qquad\qquad\qquad\qquad\;\; 
     [TH_3,F_2] =-T^2F_2, \nonumber\\
&& [T^{-1}E_1,F_1]={1\over 2}(T+T^{-1})H_1+{1\over 2}(T-T^{-1})H_2,
  \nonumber\\
&& [T^{-1}E_2,F_2]={1\over 2}(T+T^{-1})H_2+{1\over 2}(T-T^{-1})H_1, 
     \nonumber\\
&& [T^{-1}E_1,F_2] =0,\qquad\qquad\qquad\qquad\qquad\qquad\; \;
     [T^{-1}E_2,F_1] =0, \nonumber \\
  && [E_1,E_2] ={1\over 2{\sf h}}(T^2-1),\nonumber\\
  && [TF_2,TF_1]=T\biggl(F_3-{{\sf h}\over 2}H_3TH_3
                 -{{\sf h}\over 8}(T-T^{-1})\biggr) \nonumber\\
  && [TH_1,T]= {1\over 2}(T^2-1),\qquad\qquad\qquad\qquad\;\;\;\;
     [TH_1,T^{-1}]= {1\over 2}(T^{-2}-1), \nonumber\\  
  && [TH_2,T]= {1\over 2}(T^2-1),\qquad\qquad\qquad\qquad\;\;\;\;
     [TH_2,T^{-1}]= {1\over 2}(T^{-2}-1), \nonumber\\
  && [H_1,F_3]= -{T^{-1}\over 4}\biggl(TF_3+F_3T+T^{-1}F_3
                +F_3T^{-1}\biggl)-{{\sf h}\over 4}T^{-1}H_3^2
                -{{\sf h} \over 4}H_3T^{-1}H_3, \nonumber \\
  && [H_2,F_3]= -{T^{-1}\over 4}\biggl(TF_3+F_3T+T^{-1}F_3
                +F_3T^{-1}\biggl)-{{\sf h}\over 4}T^{-1}H_3^2
                -{{\sf h} \over 4}H_3T^{-1}H_3, \nonumber \\
  && [E_1,T]=[E_1,T^{-1}]=[E_2,T]=[E_2,T^{-1}]=0,\nonumber\\
  && [F_1,T]= {\sf h}TE_2, \qquad\qquad\qquad\qquad\qquad\;\;\; 
     [F_1,T^{-1}]= -{\sf h}T^{-1}E_2, \nonumber\\
  && [F_2,T]= -{\sf h}TE_1, \qquad\qquad\qquad\qquad\qquad
     [F_2,T^{-1}]={\sf h}T^{-1}E_1,\nonumber\\
  && [E_1,F_3]=-{1\over 2}\biggl(TF_2+F_2T\biggr), \qquad\qquad\;\;\;\;\;
     [E_2,F_3]={1\over 2}\biggl(TF_1+F_1T\biggr), \nonumber \\
&& [F_1,F_3]={\sf h}TF_1-{\sf h}E_2F_3+{{\sf h}^2\over 4}TE_2,\nonumber\\
  && [F_2,F_3]={\sf h}TF_2+{\sf h}E_1F_3-{{\sf h}^2\over 4}TE_1.
\end{eqnarray}
\end{prop}

Here we quoted only the final results. To obtain the realizations of
$H_1$ and $H_2$ given in (11), we, in analogy with (6), started with the 
ansatz $\sqrt{1+{\sf{h}}^2e_3^2}h_1$ and 
$\sqrt{1+{\sf{h}}^2e_3^2}h_2$ for these generators respectively. It is 
easy to see that         
\begin{eqnarray}
&& [\sqrt{1+{\sf{h}}^2e_3^2}h_1,F_3]=-{1\over 4}\biggl(TF_3+F_3T+T^{-1}F_3
   +F_3T^{-1}\biggr)\nonumber \\
&& \phantom{[\sqrt{1+{\sf{h}}^2e_3^2}h_1,F_3]=}+{{\sf{h}}^2\over 4}
\biggl(e_3(h_1-h_2)H_3 +H_3e_3(h_1-h_2)\biggr), \nonumber \\
&& [\sqrt{1+{\sf{h}}^2e_3^2}h_2,F_3]=-{1\over 4}\biggl(TF_3+F_3T+T^{-1}F_3
   +F_3T^{-1}\biggr)\nonumber \\
&& \phantom{[\sqrt{1+{\sf{h}}^2e_3^2}h_2,F_3]=}-{{\sf{h}}^2\over 4}
\biggl(e_3(h_1-h_2)H_3+H_3e_3(h_1-h_2)\biggr).
\end{eqnarray}
Then, if we add to $\sqrt{1+{\sf{h}}^2e_3^2}h_1$ and deduct from 
$\sqrt{1+{\sf{h}}^2e_3^2}h_2$ the term ${{\sf{h}}\over 2}e_3(h_1-h_2)$, 
we obtain
\begin{eqnarray}
&& [(\sqrt{1+{\sf{h}}^2e_3^2}h_1+{{\sf{h}}\over 2}e_3(h_1-h_2)),F_3]
=-{1\over 4}\biggl(TF_3+F_3T+T^{-1}F_3+F_3T^{-1}\biggr)\nonumber \\
&&\phantom{[\sqrt{1+{\sf{h}}^2e_3^2}h_1+{{\sf{h}}\over 2}e_3(h_1-h_2),F_3]=}
+{{\sf{h}}\over 4}T(h_1-h_2)H_3+{{\sf{h}}\over 4}H_3T(h_1-h_2),\nonumber\\
&&[(\sqrt{1+{\sf{h}}^2e_3^2}h_2-{{\sf{h}}\over 2}e_3(h_1-h_2)),F_3]
=-{1\over 4}\biggl(TF_3+F_3T+T^{-1}F_3+F_3T^{-1}\biggr)\nonumber \\
&&\phantom{[\sqrt{1+{\sf{h}}^2e_3^2}h_2-{{\sf{h}}\over 2}e_3(h_1-h_2),F_3]=}
-{{\sf{h}}\over 4}T(h_1-h_2)H_3-{{\sf{h}}\over 4}H_3T(h_1-h_2). 
\end{eqnarray}
These commutation relations suggest the 
realizations $H_1\sim\biggl(\sqrt{1+{\sf{h}}^2
e_3^2}h_1+{{\sf{h}} \over 2}e_3(h_1-h_2)\biggr)$ and 
$H_2\sim\biggl(\sqrt{1+
{\sf{h}}^2e_3^2}h_2-{{\sf{h}} \over 2}e_3(h_1-h_2)\biggr)$. Finally, to 
preserve the Cartan subalgebra, we are obliged to multiply both of these
expressions by $T^{-1}$. The resultant maps for  $H_1$ and $H_2$ are 
quoted in (11). The expressions of $F_1$ and $F_2$ are 
obtained in a similar way. The expressions (6), (10) and (11) may be
looked now as a particular realization of the ${\cal U}_{\sf h}(sl(3))$ 
generators. Other maps may also be considered.

\begin{prop} In terms of the Chevalley generators 
(simple roots) $\{E_1, E_2,F_1, 
F_2,H_1,H_2\}$, the algebra ${\cal U}_{{\sf{h}}}(sl(3))$ is defined as 
follows:
\begin{eqnarray}
&&T=\biggl(1+2{\sf h}[E_1,E_2]\biggr)^{1/2},\qquad\qquad\qquad\;\;
T^{-1}=\biggl(1+2{\sf h}[E_1,E_2]\biggr)^{-1/2},\nonumber\\
&&[H_1,H_2]=0, \nonumber \\
&& [H_1,E_1] =2E_1,\qquad\qquad\qquad\qquad\qquad\;\;\;
     [H_2,E_2] =2E_2,\nonumber \\
&& [H_1,E_2]=-E_2,\qquad\qquad\qquad\qquad\qquad\;\;
     [H_2,E_1]=-E_1, \nonumber\\ 
&& [H_1,F_1] =-2F_1+{\sf h}E_2(H_1+H_2), \qquad\;\;\;
     [H_2,F_2] =-2F_2-{\sf h}E_1(H_1+H_2), \nonumber\\
&& [H_1,F_2] =F_2-{\sf h}E_1(H_1+H_2), \qquad\qquad\;
     [H_2,F_1] =F_1+{\sf h}E_2(H_1+H_2), \nonumber\\
&& [T^{-1}E_1,F_1]={1\over 2}(T+T^{-1})H_1+{1\over 2}(T-T^{-1})H_2,
      \nonumber\\
&& [T^{-1}E_2,F_2]={1\over 2}(T+T^{-1})H_2+{1\over 2}(T-T^{-1})H_1, 
     \nonumber\\
&& [T^{-1}E_1,F_2] =[T^{-1}E_2,F_1] =0, \nonumber \\
&& E_1^2E_2-2E_1E_2E_1+E_2E_1^2=0,\nonumber \\
&& E_2^2E_1-2E_2E_1E_2+E_1E_2^2=0,\nonumber \\
&& (TF_1)^2TF_2-2TF_1TF_2TF_1+TF_2(TF_1)^2=0,
   \nonumber \\
&& (TF_2)^2TF_1-2TF_2TF_1TF_2+TF_1(TF_2)^2=0,
\end{eqnarray}
or, briefly
\begin{eqnarray}
&& [H_i,H_j]=0, \nonumber \\
&& [H_i,E_j] =a_{ij}E_j,\qquad\qquad\qquad
   [H_i,F_j] =-a_{ij}F_j+T^{-1}[F_j,T](H_1+H_2),\nonumber \\
&& [T^{-1}E_i,F_j]=\delta_{ij}\biggl(T^{-1}H_i+{1\over 2}
(T-T^{-1})(H_1+H_2)\biggr), \nonumber\\ 
&& (\hbox{ad}\; E_i)^{1-a_{ij}}(E_j)=0,\qquad i\neq j,\nonumber \\
&& (\hbox{ad}\; TF_i)^{1-a_{ij}}(TF_j)=0,\qquad i\neq j,
\end{eqnarray}
where $(a_{ij})_{i,j=1,2}$ is the Cartan matrix of $sl(3)$, i.e. $a_{11}=
a_{22}=2$ and $a_{12}=a_{21}=-1$. 
\end{prop}

\smallskip

{\bf 3.} We now turn to the coalgebraic structure: 

\begin{prop}
The Jordanian quantum 
algebra ${\cal U}_{\sf h}(sl(3))$ admits a Hopf structure with coproducts, 
antipodes and counits determined by
\begin{eqnarray}
&&\Delta (E_1)=E_1\otimes 1+T\otimes E_1,\nonumber \\
&& \Delta (E_2)=E_2\otimes 1+T\otimes E_2, \nonumber\\
&&\Delta (T)=T\otimes T, \qquad\qquad\qquad\qquad\qquad
\Delta (T^{-1})=T^{-1}\otimes T^{-1},\nonumber \\
&&\Delta(F_1)=F_1\otimes 1+T^{-1}\otimes F_1+{\sf h}H_3\otimes E_2\nonumber \\
&& \phantom{\Delta(F_1)}=
F_1\otimes 1+T^{-1}\otimes F_1+T(H_1+H_2)\otimes T^{-1}[F_1,T],\nonumber\\
&&\Delta(F_2)=F_2\otimes 1+T^{-1}\otimes F_2-{\sf{h}}H_3\otimes E_1 \nonumber\\
&& \phantom{\Delta(F_2)}=
F_2\otimes 1+T^{-1}\otimes F_2+T(H_1+H_2)\otimes T^{-1}[F_2,T],\nonumber\\
&&\Delta (F_3)=F_3\otimes T+T^{-1}\otimes F_3, \nonumber \\
&&\Delta (H_1)=H_1\otimes 1+1\otimes H_1
   -{1\over 2}(1-T^{-2})\otimes T^{-1}H_3\nonumber \\
&& \phantom{\Delta(H_1)}=H_1\otimes 1+1\otimes H_1
   -{1\over 2}(1-T^{-2})\otimes (H_1+H_2), \nonumber \\
&&\Delta (H_2)=H_2\otimes 1+1\otimes H_2
   -{1\over 2}(1-T^{-2})\otimes T^{-1}H_3\nonumber \\
&& \phantom{\Delta(H_2)}=H_2\otimes 1+1\otimes H_2
   -{1\over 2}(1-T^{-2})\otimes (H_1+H_2),\nonumber \\
&&\Delta (H_3)=H_3\otimes T+T^{-1}\otimes H_3,  \nonumber \\
&& S(E_1)=-T^{-1}E_1,\qquad\qquad\qquad\qquad\qquad 
   S(E_2)=-T^{-1}E_2,\nonumber\\
&& S(T)=T^{-1},\qquad\qquad\qquad\qquad\qquad\qquad\; S(T^{-1})=T,\nonumber \\
&& S(F_1)=-TF_1+{\sf h}TH_3T^{-1}E_2
=-TF_1+T^2(H_1+H_2)T^{-2}[F_1,T],\nonumber \\
&& S(F_2)=-TF_2-{\sf{h}}TH_3T^{-1}E_1=
-TF_2+T^2(H_1+H_2)T^{-2}[F_2,T],\nonumber \\
&& S(F_3)=-TF_3T^{-1},\nonumber \\
&& S(H_1)=-H_1-{1\over 2}(T-T^{-1})H_3=
-H_1-{1\over 2}(T^{2}-1)(H_1+H_2), \nonumber \\
&& S(H_2)=-H_2-{1\over 2}(T-T^{-1})H_3=
-H_2-{1\over 2}(T^{2}-1)(H_1+H_2), \nonumber \\
&& S(H_3)=-TH_3T^{-1},\nonumber\\
&&\epsilon(a)=0,\qquad\forall 
a\in\biggl\{H_1,H_2,H_3,E_1,E_2,F_1,F_2,F_3\biggr \},\nonumber \\
&& \epsilon (T)=\epsilon (T^{-1})=1.  
\end{eqnarray}
\end{prop}

All the Hopf algebra axioms can be verified by direct calculations. Let us
remark that our coproducts have simpler forms as compared to Refs.
\cite{KLM,LDO1,LDO2,LMDO}.

\begin{prop}The universal ${\cal R}_{\sf h}$-matrix has the 
following form:  
\begin{eqnarray}
&& {\cal R}_{\sf h}={\cal F}_{21}^{-1}{\cal F},  
\end{eqnarray}
where
\begin{eqnarray}
&& {\cal F}=\exp\biggl({\sf h}TH_3\otimes E_3\biggr)
\exp\biggl(2{\sf h}TE_1\otimes T^{-2}E_2\biggr). 
\end{eqnarray}
\end{prop}

The ${\cal R}$-matrix properties are verified using MAPLE. The element (18) 
coincides with the universal ${\cal R}$-matrix of the Borel subalgebra and    
gives exactly the expression (3) in the representation (fund.) $\otimes$ 
(arb.).  

\smallskip

{\bf 4.} Following Drinfeld's arguments \cite{D2}, it is possible 
to construct a twist operator $G\in {\cal U}(sl(3))^{\otimes\,2}
[[{\sf h}]]$ relating the Jordanian coalgebraic structure given
by (17) with the corresponding classical coalgebraic structure.
For an invertible map $m:{\cal U}_{\sf h}(sl(3)) \rightarrow
{\cal U}(sl(3)),\,m^{-1}:{\cal U}(sl(3)) \rightarrow {\cal U}_{\sf h}
(sl(3))$, the following relations hold:
\begin{eqnarray}
(m \otimes m) \circ \Delta \circ m^{-1} ({\cal X}) &=&
G \Delta_{0} ({\cal X}) G^{-1},\nonumber\\
m \circ S \circ m^{-1} ({\cal X}) &=& g S_{0}({\cal X}) g^{-1},
\end{eqnarray} 
where ${\cal X} \in {\cal U}(sl(3))[[{\sf h}]]$ and $( \Delta_{0},
\epsilon_{0}, S_{0})$ are the coproduct, counit and the antipode 
maps of the classical ${\cal U}(sl(3))$ algebra. The transforming 
operator $g (\in {\cal U}(sl(3))[[{\sf h}]])$ and its inverse may
be expressed as
\begin{equation}
g = \mu \circ ({\hbox {id}} \otimes S_{0}) G,\qquad
g^{-1} = \mu \circ (S_{0} \otimes {\hbox {id}}) G^{-1},
\end{equation}
where $\mu$ is the multiplication map.

For the map presented here in (6), (10) and (11), we have the 
construction
\begin{eqnarray}
&& G= 1 \otimes 1 - {\frac {1}{2}} {\sf h}{\hat r} + {\frac {1}{8}}
{\sf h}^{2}\biggl[{\hat r}^{2} + 2 (e_{3} \otimes e_{3}) 
\Delta_{0}(h_{3})\biggr]\nonumber\\
&&\phantom{G=} 
- {\frac {1}{48}} {\sf h}^{3}\biggl[{\hat r}^{3} + 6 ( e_{3}
\otimes e_{3}) \Delta_{0}(h_{3}) {\hat r} - 4 {(\Delta_{0}
(e_{3}))}^{2} {\hat r}\biggr]\nonumber\\
&& \phantom{G=}+ {\frac {1}{384}} {\sf h}^{4} \biggl[ {\hat r}^{4} 
- 16 {(\Delta_{0}(e_{3}))}^{2} {\hat r}^{2} + 12 (e_{3} \otimes
e_{3}) \Delta_{0}(h_{3}) {\hat r}^{2}+ 
12 {((e_{3} \otimes e_{3}) \Delta_{0}(h_{3}))}^{2}\nonumber\\
&&\phantom{G=} 
+ 6 {({e_{3}}^{2} \otimes 1 - 1 \otimes {e_{3}}^{2})}^{2} 
\Delta_{0}(h_{3})+ 12 {(\Delta_{0}(e_{3}))}^{2}( {e_{3}}^2 \otimes 1
+ 1 \otimes {e_{3}}^{2} ) \Delta_{0}(h_{3})\nonumber\\
&&\phantom{G=}  - 8 \Delta_{0}(e_{3}) ( {e_{3}}^{3} \otimes 1
+ 1 \otimes {e_{3}}^{3} ) \Delta _{0}(h_{3})
- 10 {(\Delta_{0}(e_{3}))}^{4} \Delta_{0}( h_{3}) \biggr] 
+ O ({\sf h}^{5}),\nonumber\\
&&g = 1 + {\sf h} e_{3} {( 1 + {\sf h}^{2} {e_{3}}^{2})}
^{1/2} + {\sf h}^{2} {e_{3}}^{2},
\end{eqnarray}
where ${\hat r} = h_{3} \otimes e_{3} - e_{3} \otimes h_{3}$. 
The above twist operators, while obeying the requirement (20) for the 
full ${\cal U}(sl(3))[[{\sf h}]]$ algebra, are, however,
generated only by the elements $(e_{3}, h_{3})$, related to
the longest root. This property accounts for 
the embedding of the ${\cal U}_{\sf h}(sl(2))$ algebra
in the higher dimensional ${\cal U}_{\sf h}(sl(3))$ algebra.
The transforming operator $g$ is obtained
in (22) in a closed form. The series expansion of the twist operator 
$G$ may be developed upto an arbitrary order in ${\sf h}$. The expansion
(22) of the twist operator $G$ in powers of ${\sf h}$ satisfies 
the cocycle condition
\begin{equation}
(1 \otimes G) ({\hbox {id}} \otimes \Delta_{0}) G
= (G \otimes 1) (\Delta_{0} \otimes {\hbox {id}}) G
\end{equation}
upto the desired order. The present discussion of the twist operator
relating to the ${\cal U}_{\sf h}(sl(3))$ algebra may be easily
extended to higher dimensional Jordanian algebras. (A systematic study of 
twists for ${\cal U}_{\sf h}(sl(2))$ can be found in [21]).    

\smallskip

{\bf 5.} Let us mention that there is a $\CC$-algebra automorphism $\phi$ 
of ${\cal U}_{\sf h}(sl(3))$ such that 
\begin{equation}
\begin{array}{lll} 
\phi(T^{\pm 1})=T^{\pm 1},\qquad \qquad & 
\phi (F_3)=F_3, \qquad \qquad &\phi(H_3)=H_3, \nonumber\\
\phi(E_1)=E_2,\qquad \qquad  
&\phi (F_1)=F_2, \qquad \qquad &\phi(H_1)=H_2, \nonumber\\
\phi(E_2)=-E_1,\qquad \qquad  
&\phi (F_2)=-F_1, \qquad \qquad &\phi(H_2)=H_1.
\end{array} 
\end{equation}
(For ${\sf h}=0$, this automorphism reduces to the classical one 
$(h_1,e_1,f_1,h_2,e_2,f_2)\longrightarrow (h_2,e_2,f_2,$ $h_1,-e_1,
-f_1)$). Also there is a second $\CC$-algebra automorphism $\varphi$ 
of ${\cal U}_{\sf h}(sl(3))$ defined as:
\begin{equation}
\begin{array}{lll} 
\varphi(T^{\pm 1})=-T^{\pm 1},\qquad \qquad  
&\varphi (F_3)=-F_3, \qquad \qquad &\varphi(H_3)=-H_3, \nonumber\\
\varphi(E_1)= E_1,\qquad \qquad  
&\varphi (F_1)= F_1, \qquad \qquad &\varphi(H_1)=H_1, \nonumber\\
\varphi(E_2)=E_2,\qquad \qquad  
&\varphi (F_2)= F_2, \qquad \qquad &\varphi(H_2)=H_2. 
\end{array}
\end{equation}

\smallskip

{\bf 6.} The expressions (6), (10) and (11) permit immediate explicit 
construction of the finite-dimensional irreducible representations 
of ${\cal U}_{\sf h}(sl(3))$. For example, the three-dimensional 
irreducible representations are spanned by
\begin{eqnarray}
   && H_1=\pmatrix{1 & 0 & {{\sf h }\over 2} \cr
                0 & 0 & 0 \cr
                0 & 0 & 0\cr}, \qquad 
   E_1=\pmatrix{0 & 1 & 0 \cr
                0 & 0 & 0 \cr
                0 & 0 & 0\cr}, \qquad 
   F_1=\pmatrix{0 & 0 & 0 \cr
                1 & 0 & -{{\sf h}\over 2} \cr
                0 & 0 & 0\cr}, \nonumber \\
 &&  H_2=\pmatrix{0 & 0 & {{\sf h}\over 2} \cr
                0 & 0 & 0 \cr
                0 & 0 & -1\cr}, \qquad 
   E_2=\pmatrix{0 & 0 & 0 \cr
                0 & 0 & 1 \cr
                0 & 0 & 0\cr},\qquad 
   F_2=\pmatrix{0 & -{{\sf h}\over 2} & 0 \cr
                0 & 0 & 0 \cr
                0 & 1 & 0\cr},\nonumber\\
 &&  H_3=\pmatrix{1 & 0 & 0 \cr
                0 & 0 & 0 \cr
                0 & 0 & -1\cr}, \qquad 
   T^{\pm 1}=\pmatrix{1 & 0 & \pm {\sf h} \cr
                      0 & 1 & 0 \cr
                      0 & 0 & 1\cr}, \qquad 
   F_3=\pmatrix{0 & 0 & 0 \cr
                0 & 0 & 0 \cr
                1 & 0 & 0 \cr}, 
\end{eqnarray}
or, by    
\begin{eqnarray}
&&   H_1=\pmatrix{1 & 0 & {{\sf h}\over 2} \cr
                0 & 0 & 0 \cr
                0 & 0 & 0\cr}, \qquad 
   E_1=\pmatrix{0 & 1 & 0 \cr
                0 & 0 & 0 \cr
                0 & 0 & 0\cr}, \qquad 
   F_1=\pmatrix{0 & 0 & 0 \cr
                1 & 0 & -{{\sf h}\over 2} \cr
                0 & 0 & 0\cr}, \nonumber \\
&&   H_2=\pmatrix{0 & 0 & {{\sf h}\over 2} \cr
                0 & 0 & 0 \cr
                0 & 0 & -1\cr}, \qquad 
 E_2=\pmatrix{0 & 0 & 0 \cr
                0 & 0 & 1 \cr
                0 & 0 & 0\cr},\qquad 
   F_2=\pmatrix{0 & -{{\sf h}\over 2} & 0 \cr
                0 & 0 & 0 \cr
                0 & 1 & 0\cr},\nonumber\\
&&   H_3=\pmatrix{-1 & 0 & 0 \cr
                0 & 0 & 0 \cr
                0 & 0 & 1\cr}, \qquad 
   T^{\pm 1}=\pmatrix{-1 & 0 & \mp {\sf h} \cr
                       0 & -1 & 0 \cr
                       0 & 0 & -1\cr}, \qquad 
   F_3=\pmatrix{0 & 0 & 0 \cr
                0 & 0 & 0 \cr
                -1 & 0 & 0 \cr}. 
\end{eqnarray}
The three-irrep. (27) is simply obtained form the irrep. (26) using 
the automorphism $\varphi$. The irrep. (27) has evidently no classical 
(${\sf h}=0$) limit.

\section{${\cal U}_{\sf h}(sl(4))$: Map and ${\cal R}_{\sf h}$-matrix} 
 
The major interest of our approch is that it can be generalized for
obtaining Jordanian quantum algebras ${\cal U}_{\sf h}(sl(N))$ of higher 
dimensions. Here we illustrate our method using ${\cal U}(sl(4))$ as an 
example. Let $h_{1}=e_{11}-e_{22}\equiv h_{12}$, $h_{2}=e_{22}-e_{33}\equiv 
h_{23}$, $h_{3}=e_{33}-e_{44}\equiv h_{34}$, $e_1\equiv e_{12}$, $e_2\equiv 
e_{23}$, $e_3\equiv e_{34}$, $f_1\equiv e_{21}$, $f_2\equiv e_{32}$ and 
$f_3\equiv e_{43}$ be the standard Chevalley generators (simple roots) of 
${\cal U}(sl(4))$. The others roots obtained by action of the Weyl group are 
denoted by $e_{13}=[e_{12},e_{23}]$, $e_{14}=[e_{13},e_{34}]$, $e_{24}=
[e_{23},e_{34}]$, $e_{31}=[e_{32},e_{21}]$, $e_{41}=[e_{43},e_{31}]$, 
$e_{42}=[e_{43},e_{32}]$, $h_{13}=h_{12}+h_{23}$, $h_{14}=h_{12}+h_{23}+
h_{34}$ and $h_{24}=h_{23}+h_{34}$. As for ${\cal U}_{\sf h}(sl(3))$, the 
Jordanian deformation arises here from the longest roots, i.e. from 
$e_{14}$, $e_{41}$ and $h_{14}$. These generators are deformed as follows:   
\begin{eqnarray}
&& T={\sf h}e_{14}+\sqrt{1+{\sf h}^2e_{14}^2}, \qquad\qquad\qquad 
T^{-1}=-{\sf h}e_{14}+\sqrt{1+{\sf h}^2e_{14}^2}, \nonumber \\
&& E_{41}=e_{41}-{{\sf h}^2\over 4}e_{14}\bigl(h_{14}^2-1\bigr),
\qquad\qquad H_{14}= \sqrt{1+{\sf h}^2e_{14}^2} h_{14}, 
\end{eqnarray}             
with the well-known coproducts
\begin{eqnarray}
&&\Delta (T)=T\otimes T, \qquad\qquad\qquad\qquad\qquad\qquad
\Delta (T^{-1})=T^{-1}\otimes T^{-1},\nonumber \\
&&\Delta (E_{41})=E_{41}\otimes T+T^{-1}\otimes E_{41}, \nonumber \\
&&\Delta (H_{14})=H_{14}\otimes T+T^{-1}\otimes H_{14}.
\end{eqnarray}

By analogy with what is happen in ${\cal U}_{\sf h}(sl(3))$ algebra, 
the subsets $\bigl\{h_{12},e_{12},e_{21},e_{24},e_{42},h_{24}=h_{23}+h_{34},
e_{14},e_{41},h_{14}=h_{12}+h_{23}+h_{34}\bigr\}$ and $\bigl\{h_{13}=h_{12}+
h_{23},e_{13},e_{31},e_{34},e_{43},h_{34},e_{14},e_{41},h_{14}
=h_{12}+h_{23}+h_{34}\bigr\}$\footnote{Each subsets forms a 
${\cal U}(sl(3))$ subalgebra in ${\cal U}(sl(4))$.} are deformed exactly as 
presented above (see (10) and (11)), i.e.       
\begin{eqnarray}
&& E_{12}=\sqrt{{\sf h}e_{14}+\sqrt{1+{\sf h}^2e_{14}^2}}e_{12}=
T^{1/2}e_{12},\nonumber \\
&& E_{24}=\sqrt{{\sf h}e_{14}+\sqrt{1+{\sf h}^2e_{14}^2}}e_{24}=
T^{1/2}e_{24}, \nonumber\\
&& E_{21}=\sqrt{-{\sf h}e_{14}+\sqrt{1+{\sf h}^2e_{14}^2}}e_{21}
+{{\sf h}\over 2}\sqrt{{\sf h}e_{14}+\sqrt{1+{\sf h}^2e_{14}^2}}
e_{24}h_{14}=T^{-1/2}\biggl(e_{21}+{{\sf h}\over 2}T
e_{24}h_{14} \biggr),\nonumber \\
&& E_{42}=\sqrt{-{\sf h}e_{14}+\sqrt{1+{\sf h}^2e_{14}^2}}e_{42}-
{{\sf h}\over 2}\sqrt{{\sf h}e_{14}+\sqrt{1+{\sf h}^2e_{14}^2}}
e_{12}h_{14}=T^{-1/2}\biggl(e_{42}-{{\sf h}\over 2}T
e_{12}h_{14} \biggr), \nonumber \\
&& H_{12}=\biggl(-{\sf h}e_{14}+\sqrt{1+{\sf h}^2e_{14}^2}\biggr)
    \biggl(\sqrt{1+{\sf h}^2e_{14}^2}h_{12}
    +{{\sf h}\over 2}e_{14}(h_{12}-h_{24})\biggr)=
h_{12}-{{\sf h}\over 2}e_{14}T^{-1}h_{14}, \nonumber\\
&& H_{24}=\biggl(-{\sf h}e_{14}+\sqrt{1+{\sf h}^2e_{14}^2}\biggr)
 \biggl(\sqrt{1+{\sf h}^2e_{14}^2}h_{24}
    -{{\sf h}\over 2}e_{14}(h_{12}-h_{24})\biggr)=
h_{24}-{{\sf h}\over 2}e_{14}T^{-1}h_{14}
\end{eqnarray}
and
\begin{eqnarray}
&& E_{13}=\sqrt{{\sf h}e_{14}+\sqrt{1+{\sf h}^2e_{14}^2}}e_{13}
=T^{1/2}e_{13},\nonumber \\
&& E_{34}=\sqrt{{\sf h}e_{14}+\sqrt{1+{\sf h}^2e_{14}^2}}e_{34}
=T^{1/2}e_{34}, \nonumber\\
&& E_{31}=\sqrt{-{\sf h}e_{14}+\sqrt{1+{\sf h}^2e_{14}^2}}e_{31}+
{{\sf h}\over 2}\sqrt{{\sf h}e_{14}+\sqrt{1+{\sf h}^2e_{14}^2}}
e_{34}h_{14}=T^{-1/2}\biggl(e_{31}+{{\sf h}\over 2}e_{34}h_{14}\biggr),
\nonumber \\
&& E_{43}=\sqrt{-{\sf h}e_{14}+\sqrt{1+{\sf h}^2e_{14}^2}}e_{43}-{{\sf h}
\over 2}\sqrt{{\sf h}e_{14}+\sqrt{1+{\sf h}^2e_{14}^2}}e_{13}h_{14}
=T^{-1/2}\biggl(e_{43}-{{\sf h}\over 2}e_{13}h_{14}\biggr), \nonumber \\
&& H_{13}=\biggl(-{\sf h}e_{14}+\sqrt{1+{\sf h}^2e_{14}^2}\biggr)
\biggl(\sqrt{1+{\sf h}^2e_{14}^2}h_{13}+{{\sf h}\over 2}e_{14}
(h_{13}-h_{34})\biggr)=h_{13}-{{\sf h}\over 2}e_{14}T^{-1}h_{14}, \nonumber\\
&& H_{34}=\biggl(-{\sf h}e_{14}+\sqrt{1+{\sf h}^2e_{14}^2}\biggr)
\biggl(\sqrt{1+{\sf h}^2e_{14}^2}h_{34}-{{\sf h}\over 2}e_{14}
(h_{13}-h_{34})\biggr)=h_{34}-{{\sf h}\over 2}e_{14}T^{-1}h_{14}.
\end{eqnarray}

The elements $E_{23}$, $E_{32}$ and $H_{23}$ are obtained after 
analyzing the commutators $[E_{24},E_{43}]$ and $[E_{34},E_{42}]$. It is 
simple to see that these elements remain undeformed, i.e.
\begin{eqnarray}
&& E_{23}=e_{23},\qquad\qquad
E_{32}=e_{32},\qquad\qquad H_{23}=h_{23}.  
\end{eqnarray}
It is now easy to verify that    
\begin{eqnarray}
&& H_{23}+H_{34}=H_{24}, \qquad [E_{12},E_{23}]=E_{13}, \qquad
   [E_{32},E_{21}]=E_{31}, \nonumber \\
&& H_{12}+H_{23}=H_{13}, \qquad [E_{23},E_{34}]=E_{24}, \qquad
   [E_{43},E_{32}]=E_{42}. 
\end{eqnarray}

\begin{prop}The generating elements $H_1\equiv H_{12}$, $H_2\equiv H_{23}$,
$H_3\equiv H_{34}$, $E_1\equiv E_{12}$, $E_2\equiv E_{23}$, $E_3\equiv E_{34}$,
$F_1\equiv E_{21}$, $F_2\equiv E_{32}$, $F_3\equiv E_{43}$ of the Jordanian 
quantum algebra ${\cal U}_{\sf h}(sl(4))$ obey the following commutations 
rules:    
\begin{eqnarray}
&&T=\biggl(1+2{\sf h}[E_1,[E_2,E_3]]\biggr)^{1/2},\qquad\qquad\qquad\;\;
T^{-1}=\biggl(1+2{\sf h}[E_1,[E_2,E_3]]\biggr)^{-1/2},\nonumber\\
&&[H_1,H_2]=[H_1,H_3]=[H_2,H_3]=0, \nonumber \\
&&[H_1,E_1] =2E_1, \qquad\qquad [H_1,E_2] =-E_2,\qquad\qquad
[H_1,E_3] =0,\nonumber\\
&&[H_2,E_1]=-E_1,\qquad\qquad [H_2,E_2] =2E_2,\qquad\qquad 
[H_2,E_3]=-E_3,\nonumber\\
&&[H_3,E_1] =0,\qquad \qquad [H_3,E_2]=-E_2,\qquad\qquad
[H_3,E_3] =2E_3, \nonumber\\
&&[H_1,F_1] =-2F_1+T^{-1}[F_1,T](H_1+H_2+H_3), \qquad
[H_1,F_2] =F_2,\nonumber \\
&&[H_1,F_3] =T^{-1}[F_3,T](H_1+H_2+H_3),\nonumber\\
&&[H_2,F_1]=F_1,\qquad [H_2,F_2] =-2F_2,\qquad [H_2,F_3]=F_3,\nonumber\\
&&[H_3,F_1] =T^{-1}[F_1,T](H_1+H_2+H_3),\qquad [H_3,F_2]=F_2,\nonumber \\
&&[H_3,F_3] =-2F_3+T^{-1}[F_3,T](H_1+H_2+H_3), \nonumber\\
&&[T^{-1}E_1,F_1]=T^{-1}H_1+{1\over 2}(T-T^{-1})(H_1+H_2+H_3),\nonumber\\
&& [E_2,F_2]=H_2,\nonumber\\
&& [T^{-1}E_3,F_3]=T^{-1}H_3+{1\over 2}(T-T^{-1})(H_1+H_2+H_3),\nonumber\\
&& [T^{-1}E_1,F_2] =[T^{-1}E_1,F_3] =0, \nonumber \\
&& [E_2,F_1] =[E_2,F_3] =0, \nonumber \\
&& [T^{-1}E_3,F_1] =[T^{-1}E_3,F_2] =0, \nonumber \\
&& [E_1,E_3]=[TF_1,TF_3]=0,\nonumber \\
&& E_1^2E_2-2E_1E_2E_1+E_2E_1^2=0,\qquad\qquad\qquad\;
E_1E_2^2-2E_2E_1E_2+E_2^2E_1=0,\nonumber \\
&& E_2^2E_3-2E_2E_3E_2+E_3E_2^2=0,\qquad\qquad\qquad\;
E_2E_3^2-2E_3E_2E_3+E_3^2E_2=0,\nonumber \\
&& (TF_1)^2F_2-2TF_1F_2TF_1+F_2(TF_1)^2=0,\qquad
TF_1F_2^2-2F_2TF_1F_2+F_2^2TF_1=0,\nonumber \\
&& (TF_3)^2F_2-2TF_3F_2TF_3+F_2(TF_3)^2=0,\qquad
F_2^2TF_3-2F_2TF_3F_2+TF_3F_2^2=0,
\end{eqnarray}
or, briefly,
\begin{eqnarray}
&&[H_i,H_j]=0, \nonumber \\
&&[H_i,E_j] =a_{ij}E_j, \nonumber\\
&&[H_i,F_j] =-a_{ij}F_j+(\delta_{i1}+\delta_{i3})
T^{-1}[F_j,T](H_1+H_2+H_3),\nonumber \\
&&[T^{-(\delta_{i1}+\delta_{i3})}E_i,F_j]=\delta_{ij}\biggl(T^{-(\delta_{i1}+
\delta_{i3})}H_i+{(\delta_{i1}+\delta_{i3})\over 2}
(T-T^{-1})(H_1+H_2+H_3)\biggr),\nonumber\\
&& [E_i,E_j]=[T^{(\delta_{i1}+\delta_{i3})}F_i,
T^{(\delta_{j1}+\delta_{j3})}F_j]=0,\qquad\qquad |i-j|>1,\nonumber \\
&& (\hbox{ad}\; E_i)^{1-a_{ij}}(E_j)=0,\qquad (i\neq j),\nonumber \\
&&(\hbox{ad}\; T^{(\delta_{i1}+\delta_{i3})}F_i)^{1-a_{ij}}
(T^{(\delta_{j1}+\delta_{j3})}F_j)=0,\qquad (i\neq j),
\end{eqnarray}
where $(a_{ij})_{i,j=1,2,3}$ is the Cartan matrix of $sl(4)$. 
\end{prop}

\begin{prop} The non-cocommutative coproduct structure of ${\cal U}_{\sf h}
(sl(4))$ reads:
\begin{eqnarray}
&&\Delta (E_1)=E_1\otimes 1+T\otimes E_1,\nonumber \\
&&\Delta (E_2)=E_2\otimes 1+1\otimes E_2,\nonumber \\
&&\Delta (E_3)=E_3\otimes 1+T\otimes E_3, \nonumber \\
&&\Delta(F_1)=F_1\otimes 1+T^{-1}\otimes F_1+
(H_1+H_2+H_3)\otimes T^{-1}[F_1,T],\nonumber\\
&&\Delta(F_2)=F_2\otimes 1+T^{-1}\otimes F_2,\nonumber\\
&&\Delta(F_3)=F_3\otimes 1+T^{-1}\otimes F_3+
(H_1+H_2+H_3)\otimes T^{-1}[F_3,T],\nonumber\\
&&\Delta (H_1)=H_1\otimes 1+1\otimes H_1-{1\over 2}(1-T^{-2})
\otimes (H_1+H_2+H_3), \nonumber \\
&&\Delta (H_2)=H_2\otimes 1+1\otimes H_2,\nonumber \\
&&\Delta (H_3)=H_3\otimes 1+1\otimes H_3-{1\over 2}(1-T^{-2})
\otimes (H_1+H_2+H_3).
\end{eqnarray}
\end{prop}

In the (fund.) $\otimes$ (arb.) representation, the $R_{\sf h}=(
\pi_{(fund.)}\otimes\pi_{(arb.)}){\cal R}_{\sf h}$ take 
the following simple form:
\begin{eqnarray}
&& R_{\sf h}=\pmatrix{T 
&2{\sf h}T^{-1/2}e_{24} & 2{\sf h}T^{-1/2} e_{34}
&-{{\sf h}\over 2}(T+T^{-1})\bigl(h_1+h_2+h_3\bigr)
+{{\sf h}\over 2}\bigl(T-T^{-1}\bigr)\cr
0 & I & 0 & -2{\sf h}T^{1/2}e_{12}  \cr
0 & 0 & I & -2{\sf h}T^{1/2}e_{13}  \cr
0 & 0 & 0 & T^{-1}  \cr}.
\end{eqnarray}

\begin{prop} The universal ${\cal R}_{\sf h}$-matrix for ${\cal U}_{\sf h}
(sl(4))$ may be cast in the form:             
\begin{eqnarray}
&& {\cal R}_{\sf h}={\cal F}_{21}^{-1}{\cal F}, 
\end{eqnarray}
where
\begin{eqnarray}
&&{\cal F}=\exp\biggl({\sf h}TH_{14}\otimes E_{14}\biggr)
\exp\biggl(2{\sf h}TE_{34}\otimes T^{-2}E_{13}+
2{\sf h}TE_{24}\otimes T^{-2}E_{12}\biggr), \\
&& E_{14}={\sf h}^{-1}\ln T= {\sf h}^{-1}\arcsinh{\sf h}e_{14}.
\end{eqnarray}
\end{prop}

The ${\cal R}_{\sf h}$-matrix (38) coincides with the universal 
${\cal R}$-matrix of the Borel subalgebra. Let us just note that the tensor 
elements $TE_{34}\otimes T^{-2}E_{13}$ and $TE_{24}\otimes T^{-2}E_{12}$ 
commute.

\smallskip

\section{${\cal U}_{\sf h}(sl(N))$: Generalization}

The ${\cal U}_{\sf h}(sl(5))$ algebra is derived in a similar 
way: The elements $E_2$, $E_3$, $F_2$, $F_3$, $H_2$, $H_3$ are 
not affected by the nonstandard quantization. From these above studies, It 
is easy to see that: 

\begin{prop}
The Jordanian quantization deform ${\cal U}_{\sf h}(sl(N))$'s Chevalley 
generators as follows: 
\begin{eqnarray}
&& T={\sf h}[e_1,[e_2,\cdots,[e_{N-2},e_{N-1}]\cdots]]+
\sqrt{1+{\sf h}^2([e_1,[e_2,\cdots,[e_{N-2},e_{N-1}]\cdots]])^2}, \nonumber \\
&& T^{-1}=-{\sf h}[e_1,[e_2,\cdots,[e_{N-2},e_{N-1}]\cdots]]+
\sqrt{1+{\sf h}^2([e_1,[e_2,\cdots,[e_{N-2},e_{N-1}]\cdots]])^2}, \nonumber \\ 
&&E_i=T^{(\delta_{i1}+\delta_{i,N-1})/2}e_i,\nonumber\\
&& F_i=T^{-(\delta_{i1}+\delta_{i,N-1})/2}\biggl(
f_i+{{\sf h}\over 2}T[f_i,[e_1,[e_2,\cdots,[e_{N-2},e_{N-1}]\cdots]]](h_1+\cdots +
h_{N-1})\biggr)\nonumber\\
&& H_i=h_i-{(\delta_{i1}+\delta_{i,N-1}){\sf h}\over 2}
[e_1,[e_2,\cdots,[e_{N-2},e_{N-1}]\cdots]]T^{-1}(h_1+\cdots +
h_{N-1}) 
\end{eqnarray}
($i=1,\cdots,\;N-1$) and they satisfy the commutation relations
\begin{eqnarray}
&& [H_i,H_j]=0, \nonumber \\
&& [H_i,E_j] =a_{ij}E_j,\nonumber \\
&& [H_i,F_j] =-a_{ij}F_j+(\delta_{i1}+\delta_{i,N-1})T^{-1}[F_j,T]
(H_1+\cdots + H_{N-1}), \nonumber\\
&& [T^{-(\delta_{i1}+\delta_{i,N-1})}E_i,F_j]=
\delta_{ij}\biggl(T^{-(\delta_{i1}+\delta_{i,N-1})}H_i
+{(\delta_{i1}+\delta_{i,N-1})
\over 2}(T-T^{-1})(H_1+\cdots +H_{N-1})\biggr),\nonumber\\
&& [E_i,E_j] =0,\qquad\qquad \qquad\qquad |i-j|>1,\nonumber \\ 
&& [T^{(\delta_{i1}+\delta_{i,N-1})}F_i,T^{(\delta_{j1}+\delta_{j,N-1})}F_j]=0,
\qquad\qquad |i-j|>1,\nonumber \\ 
&& (\hbox{ad} \;E_i)^{1-a_{ij}}(E_j)=0,\qquad (i\neq j),\nonumber \\
&&(\hbox{ad} \; T^{(\delta_{i1}+\delta_{i,N-1})}F_i)^{1-a_{ij}}
(T^{(\delta_{j1}+\delta_{j,N-1})}F_j)=0,\qquad (i\neq j),
\end{eqnarray} 
where $(a_{ij})_{i,j=1,\cdots , N}$ is the Cartan matrix of $sl(N)$, 
i.e. $a_{ii}=2$, $a_{i,i\pm 1}=-1$ and $a_{ij}=0$ for $|i-j|>1$.    
\end{prop}

The algebra (42) is called the {\it Jordanian quantum algebra ${\cal U}_{\sf h}
(sl(N))$}. 
The expressions (41) may be regarded as a particular nonlinear realization of 
the ${\cal U}_{\sf h}(sl(N))$ generators.   

\begin{prop} The Jordanian algebra ${\cal U}_{\sf h}(sl(N))$ (42) 
admits the following coalgebra structure:  
\begin{eqnarray}
&&\Delta (E_i)=E_i\otimes 1+T^{(\delta_{i1}+\delta_{i,N-1})}
\otimes E_i, \nonumber \\
&&\Delta (F_i)=F_i\otimes 1+T^{-(\delta_{i1}+\delta_{i,N-1})} 
\otimes F_i+T(H_1+\cdots +H_{N-1})\otimes T^{-1}[F_i,T],\nonumber \\
&&\Delta (H_i)=H_i\otimes 1+1\otimes H_i-
{(\delta_{i1}+\delta_{i,N-1})\over 2}(1-T^{-2})\otimes (H_1+\cdots +H_{N-1}), 
\nonumber\\
&& S(E_i)=-T^{-(\delta_{i1}+\delta_{i,N-1})}E_i, \nonumber \\
&& S(F_i)=-T^{(\delta_{i1}+\delta_{i,N-1})}F_i+T^2(H_1+\cdots +H_{N-1})
T^{-2}[F_i,T], \nonumber \\
&& S(H_i)=-H_i+{(\delta_{i1}+\delta_{i,N-1})\over 2}(1-T^2)
(H_1+\cdots +H_{N-1}), \nonumber \\
&& \epsilon(E_i)= \epsilon(F_i)=\epsilon(H_i)=0.
\end{eqnarray}
\end{prop}

 \begin{prop} The ${\cal R}_{\sf h}$-matrix of 
${\cal U}_{\sf h}(sl(N))$ has the following general form:
\begin{eqnarray}
&& {\cal R}_{\sf h}={\cal F}_{21}^{-1}{\cal F}, 
\end{eqnarray}
where
 \begin{eqnarray} 
&& {\cal F}=\exp\biggl({\sf h}TH_{1N}\otimes E_{1N}\biggr)
            \exp \biggl(\sum_{k=2}^{N-1}2{\sf h}TE_{kN}\otimes 
              T^{-2}E_{1k}\biggr), \\
&& H_{1N}=T(H_1+\cdots H_{N-1}), \\
&& E_{1N}={\sf h}^{-1}\ln T= {\sf h}^{-1}\arcsinh{\sf h}e_{1N}, \\
&& E_{kN}=[E_k,[\cdots,[E_{N-2},E_{N-1}]]], \qquad k=2,\cdots,N-2, \\
&& E_{N-1,N}=E_{N-1}, \\
&& E_{12}=E_1,\\
&& E_{1k}=[E_1,[\cdots,[E_{k-2},E_{k-1}]]], \qquad k=3,\cdots,N-1
\end{eqnarray}
and may be obtained from the ${\cal R}_q$-matrix associated to 
${\cal U}_q(sl(N))$ via the contraction procedure discussed above, i.e.
\begin{eqnarray} 
{\cal R}_{\sf h}=\lim_{q\rightarrow 1}
\biggl[E_q\biggl({{\sf h}{\hat e}_{1N}\over q-1}\biggr)\otimes 
E_q\biggl({{\sf h}{\hat e}_{1N}\over q-1}\biggr)\biggr]^{-1}
{\cal R}_q\biggl[E_q\biggl({{\sf h}{\hat e}_{1N} \over q-1}\biggr)
\otimes E_q\biggl({{\sf h}{\hat e}_{1N}\over q-1}\biggr)\biggr].
\end{eqnarray} 
\end{prop}

It is interesting to note that, via the nonlinear map (41), the 
${\sf h}$-deformed generators $(E_i,F_i,H_i)$ may be also equipped with an 
induced co-commutative coproduct. Similarly, the undeformed generators
$(e_i,f_i,h_i)$, via the inverse map, may be viewed as elements of the 
${\cal U}_{\sf h}(sl(N))$ algebra; and, thus, may be endowed with an 
induced noncommutative coproduct.

\vskip 1cm 

\noindent $\underline{\hbox{\bf Acknowledgments}}$: One of us (BA) wants 
to thank Professor 
Peter Forgacs for a kind invitation to the University of Tours, where parts 
of this work was done. He is also grateful to the members of the group 
for their kind hospitality.  

\smallskip
\smallskip   

\noindent $\underline{\;\;\;\;\;\;\;\;\;\;\;\;\;\;\;\;\;\;\;\;\;\;\;\;\;\;
\;\;\;\;\;\;\;\;\;\;\;\;\;\;\;\;\;\;\;\;\;\;\;\;\;\;\;\;\;\;\;\;\;\;\;\;\;
\;\;\;\;\;\;\;\;\;\;\;\;\;\;\;\;\;\;\;\;\;\;\;\;\;\;\;\;\;\;\;\;\;\;\;\;\;
\;\;\;\;\;\;\;\;\;\;\;\;\;\;\;\;\;\;\;\;\;\;\;\;\;\;\;\;\;\;\;\;\;\;\;\;\;
\;\;\;\;\;\;}$


\begin{thebibliography}{99}


\bibitem {D1} V. G. Drinfeld, {\it Quantum Groups,} Proc.
  Int. Congress of Mathematicians, Berkeley, California, Vol.
  {\bf 1}, Academic Press, New York (1986), 798.

\bibitem{J} M. Jimbo, {\it A $Q$ Difference analog of ${\cal U}(G)$ and
the Yang-Baxter equation,} 
Lett. Math. Phys. {\bf 10} (1985) 63-60.

\bibitem {DMMZ} E. E. Demidov, Yu. I. Manin, E. E. Mukhin, D. Z. Zhdanovich,
{\it Nonstandard quantum deformation of $GL(N)$ and constant solutions of 
the Yang-Baxter equation,} Prog. Theor. Phys. Suppl. 102 (1990) 203-218.  

\bibitem {O} Ch. Ohn, {\it A $\star$-product on $SL(2)$ and the 
corresponding nonstandard quantum group $GL(2)$},    
Lett. Math. Phys. {\bf 25} (1992) 85-88.

\bibitem{D2} V. G. Drinfeld, Leningrad Math. J. {\bf 1} (1990) 1419. 

\bibitem{Og} O. V. Ogievetsky, Suppl. Rendiconti Cir. Math. Palmero, Serie 
II {\bf 37}, (1993) 4569. 

\bibitem{GGS} M. Gerstenhaber, A. Giaquinto and S. D. Schak, Israel Math. 
Conf. Proc. {\bf 7} (1993) 45.      

\bibitem {KLM} P. P. Kulish, V.D. Lyakhovsky and A.I. Mudrov, 
{\it Extended Jordanian twists for Lie algebras}, J. Math. Phys. {\bf 40}
(1999) 4569.

\bibitem {LDO1} V. D. Lyakhovsky and M. A. Del Olmo, 
{\it Peripheric Extended Twists}, math.QA/9811153. 

\bibitem{LDO2} V. D. Lyakhovsky and M. A. Del Olmo, 
{\it Extended and Reshetikhin Twists for $sl(3)$}, math.QA/9903065.

\bibitem{LMDO} V. D. Lyakhovsky, A. M. Mirolubov and M. A. Del Olmo, 
{\it Quantum Jordanian Twists}, math.QA/9811153.
 
\bibitem{AKL} D. N. Ananikian, P. P. Kulish and V.D.Lyakhovsky, 
{\it Chains of twists for symplectic Lie algebras}, math.QA/0010312.

\bibitem{KLDO} P. P. Kulish, V.D. Lyakhovsky and M. A. Del Olmo, 
{\it Chains of twists for classical Lie algebras}, math.QA/9908061.

\bibitem{KL} P. P. Kulish and V.D. Lyakhovsky, 
{\it Jordanian twists on deformed carrier subspaces}, math.QA/0007182.

 \bibitem{KLS} P. P. Kulish, V.D. Lyakhovsky and A. Stolin, 
{\it Full chains of twists for orthogonal algebras}, math.QA/0007182.

\bibitem{LSK} V.D. Lyakhovsky, A. Stolin and P. P. Kulish, 
{\it Chains of Frobenius subalgebras of $so(M)$ and the corresponding twists},
math.QA/0010147.
 
\bibitem{CK} E. Celeghini and P. P. Kulish, {\it Twist Deformation of the rank
one Lie Superalgebra}, J. Phys. A: {31}, (1998) L79. 

\bibitem{K} P. P. Kulish, {\it Super-Jordanian deformation of the 
orthosymplectic Lie superalgebras}, math.QA/9806104. 

\bibitem{ACC1} B. Abdesselam, A. Chakrabarti and R. Chakrabarti, 
{\it Irreducible representations of Jordanian quantum algebra 
${\cal U}_h(sl(2))$ via a nonlinear map}, 
Mod. Phys. Lett. A {\bf 36} (1996) 2883.

\bibitem{ACC2} B. Abdesselam, A. Chakrabarti and R. Chakrabarti,
{\it Towards a general construction of nonstandard $R_h$-matrices as 
contraction limits of $R_q$-matrices: The ${\cal U}_h(sl(N))$ algebra case}, 
Mod. Phys. Lett. A {\bf 10} (1998) 779. 

\bibitem{ACCS} B. Abdesselam, A. Chakrabarti, R. Chakrabarti and J. Segar,
{\it Maps and twists relating ${\cal U}(sl(2))$ and the nonstandard 
${\cal U}_h(sl(2))$: unified construction}, 
Mod. Phys. Lett. A {\bf 12} (1999) 765-777. 

\bibitem{Ali} M. Alishahiha, J. Phys. A {\bf 28} (1995) 6187
\bibitem{M} S. Majid, {\it Foundations of Quantum Group Theory} (Cambridge 
Univ. Press, 1995).  

\end{thebibliography}
\end{document}